# A MICROSCOPIC MODEL FOR STEFAN'S MELTING AND FREEZING PROBLEM

By Claudio Landim and Glauco Valle

*IMPA and CNRS, and Ecole Polytechnique Fédérale de Lausanne*

We study a class of one-dimensional interacting particle systems with random boundaries as a microscopic model for Stefan's melting and freezing problem. We prove that under diffusive rescaling these particle systems exhibit a hydrodynamic behavior described by the solution of a Cauchy–Stefan problem.

**1. Introduction.** In this work we return to the classical Stefan's freezing on the ground model [8]. It could be described in the following way: Consider the real line occupied by a heat-conducting material (heat is transmitted only by conduction). This material is initially almost everywhere characterized by a bounded and measurable temperature function $T:\mathbb{R} \to \mathbb{R}$. According to the temperature the material could be in one of two phases, a liquid phase for positive temperatures and a solid phase for negative temperatures. The temperature $T=0$ is that of crystallization at which both phases may occur. The problem consists in determining the temporal evolution of the temperature profile.

We consider this problem under more restrictive conditions. Suppose that at initial time the liquid phase fills the domain $u>0$ at positive temperatures and the solid phase fills the domain $u<0$ at negative temperatures. Denote by $\rho_{-1}^0 : \mathbb{R}_- \to \mathbb{R}_-$ and $\rho_1^0 : \mathbb{R}_+ \to \mathbb{R}_+$ the initial temperature profile. We are able to determine a function $B = B(t)$ describing the time evolution of the boundary between the two phases and their temperature functions, respectively $\rho_{-1}(t,u)$ and $\rho_1(t,u)$ for the solid and liquid phases. It is well known that these functions satisfy a Cauchy–Stefan problem:

$$\partial_t \rho_{-1} = a_{-1}\partial_{uu}\rho_{-1}, \qquad \partial_t \rho_1 = a_1 \partial_{uu}\rho_1,$$
(1.1) $$\dot{B}(t) = k\{a_{-1}\partial_u \rho_{-1}(t,B(t)) - a_1 \partial_u \rho_1(t,B(t))\}, \qquad B(0)=0,$$

$$\rho_i(t, B(t)) = 0, \qquad \rho_i(0, \cdot) = \rho_i^0(\cdot),$$

where $\rho_{-1}^0 : \mathbb{R}_- \to \mathbb{R}_-$ and $\rho_1^0 : \mathbb{R}_+ \to \mathbb{R}_+$ are bounded measurable functions, $a_{-1} \geq 0$ and $a_1 > 0$ are the coefficients of heat conduction of the material with respect to the solid and liquid phases and $k > 0$ is a scaling factor for the temperature.

In this paper we present a microscopic model for Stefan's equation through an appropriated interacting particle system and scaling limit techniques. Such sorts of descriptions have been proposed previously by Chayes and Swindle [2] in the case of finite domains with $a_{-1} = 0$. Rezakhanlou [6] and Bertini, Buttà and Rüdiger [1] also derive Stefan's equation as hydrodynamic limit of interacting particle systems.

We shall denote by $\mathbb{Z}$, $\mathbb{N}$ and $\mathbb{Z}_-$ the sets of integers, positive integers and nonpositive integers respectively, and by $\mathbb{R}_-$ and $\mathbb{R}_+$ the sets of nonpositive and nonnegative real numbers.

For the informal description of the microscopic model, consider the one-dimensional lattice $\mathbb{Z}$ with each site being occupied by a molecular agglomerate of type $-1$ for the material in the solid state and of type 1 for the material in the liquid state. According to its internal energy, each agglomerate is classified by a heat unit of 0 or 1. An interaction between neighboring sites occurs independently in the following way: If the particles are of the same type, then their heat units are interchanged after a mean $a_{-1}$ exponential time for particles type $-1$ and after a mean $a_1$ exponential time for particles of type 1. If the particles are of distinct type and their heat units are also distinct, at rate 1 the heat unit of the agglomerate whose heat value was 1 drops to 0 and simultaneously the other agglomerate changes type. If the particles are of distinct type and their heat units are equal to 1, both heat units drop to 0 after a mean 1 exponential time. Moreover, we start with configurations such that the agglomerates are of type $-1$ if they are at the left of the origin; otherwise they are of type 1.

We will show that this system has a hydrodynamic behavior under diffusive scaling described by the solution of a Cauchy–Stefan problem of type (1.1) with scaling factor $k = 1$, where the temperature is the macroscopic heat density profile. The general case with an arbitrary $k$ can be obtained from the previous one rescaling the temperature by $k^{-1}$. Here the diffusive scaling is expected since the hydrodynamic behavior of the simple symmetric exclusion process is described in this scale. Actually, our model could be described as a coupling between two one-dimensional nearest-neighbor simple symmetric exclusion processes in the semi-infinite lattice. To make this identification, consider each agglomerate of solid (resp. liquid) phase as a site in the space of the $\mathbb{Z}_-$ (resp. $\mathbb{N}$) in such a way that this association preserves the order. At each site whose associated agglomerate has heat unit equal to 1 we put a particle. In each of the lattices $\mathbb{Z}_-$, $\mathbb{N}$ particles evolve as in a nearest-neighbor one-dimensional simple symmetric exclusion process with jump rates $a_{-1}$ and $a_1$, respectively. Superposed to this dynamics, a particle at the boundary of one of the lattices waits a mean 1 exponential time and attempts to leave the system. If no particle occupies the boundary of the other lattice, this particle vanishes triggering a translation of the whole system to the right or left depending whether the particle was occupying a site in $\mathbb{Z}_-$ or in $\mathbb{N}$. If the boundary of the other lattice is occupied, then both particles leave the system.

From the technical point of view, the main difficulty of this problem lies in the fact that no entropy argument can be used due to the annihilation mechanism. Indeed, if one fixes the boundary to be at the origin, the unique invariant measure is degenerate, and even when estimating the relative entropy with respect to a nonstationary state as in [3, 5], the translations introduce expressions too large to be estimated by the sole entropy. Coupling is therefore the unique available tool in this context.

**2. Main result.** Let $\Gamma = \{(-1,0),(-1,1),(1,0),(1,1)\}^{\mathbb{Z}}$ be the configuration space and denote a typical configuration by $(\sigma,\eta) = \{(\sigma(x),\eta(x))\}_{x\in\mathbb{Z}} \in \Gamma$. Fix $a_{-1} \geq 0$, $a_1 > 0$ and define the generators $\mathcal{G}_0$, $\mathcal{G}_1$ by

$$(\mathcal{G}_0 f)(\sigma,\eta) = \sum_{x\in\mathbb{Z}} \mathbf{1}\{\sigma(x) = \sigma(x+1)\} a_{\sigma(x)} \{f(\sigma,\eta^{x,x+1}) - f(\sigma,\eta)\},$$

$$(\mathcal{G}_1 f)(\sigma,\eta) = \sum_{x\in\mathbb{Z}} \mathbf{1}\{\sigma(x) \neq \sigma(x+1)\} \{f(T^{x,x+1}(\sigma,\eta)) - f(\sigma,\eta)\}.$$

In these formulas, $f$ stands for a cylinder function $f:\Gamma \to \mathbb{R}$, $\eta^{x,y}$ for the configuration $\eta$ with spins at $x$, $y$ interchanged,

(2.1) $$\eta^{x,x+1}(z) = \begin{cases} \eta(x+1), & \text{if } z = x, \\ \eta(x), & \text{if } z = x+1, \\ \eta(z), & \text{otherwise,} \end{cases}$$

and $T^{x,x+1}(\sigma,\eta) = (\tilde{\sigma},\tilde{\eta})$ for the configuration defined as follows:

$$\tilde{\sigma}(z) = \begin{cases} \sigma(x), & \text{if } z = x+1, \eta(x+1) = 0, \eta(x) = 1, \\ \sigma(x+1), & \text{if } z = x, \eta(x) = 0, \eta(x+1) = 1, \\ \sigma(z), & \text{otherwise,} \end{cases}$$

$$\tilde{\eta}(z) = \begin{cases} 0, & \text{if } z = x, x+1, \\ \eta(z), & \text{otherwise.} \end{cases}$$

Hence, two neighboring particles of different type annihilate each other at rate 1: $(-1,1),(1,1) \to (-1,0),(1,0)$; and a particle sitting next to a vacant site of different type dies transforming the neighboring site into a site of its type at rate 1: $(-1,0),(1,1) \to (1,0),(1,0)$.

Denote by $(\sigma_t,\eta_t)$ the Markov process associated to the generator $\mathcal{G} = \mathcal{G}_0 + \mathcal{G}_1$ speeded up by $N^2$. The goal of this article is to show that its macroscopic behavior is described by solutions of Stefan's equations.

*The Cauchy–Stefan problem.* Let $a(u) = a_{-1}\mathbf{1}\{u \le 0\} + a_1\mathbf{1}\{u > 0\}$, $A(u) = ua(u)$. Fix a bounded measurable function $\rho_0 \colon \mathbb{R} \to \mathbb{R}$ and consider the Stefan problem

$$\partial_t \rho = a(\rho)\Delta\rho,$$
(2.2) $$\dot{u}_i(t) = b_i\{\partial_u A(\rho(t, u_i(t)-)) - \partial_u A(\rho(t, u_i(t)+))\},$$
$$\rho(0, \cdot) = \rho_0(\cdot).$$

Here $\Delta$ stands for the Laplacian, $u_i(t)$ for the curves at which $\rho(t, u_i(t)) = 0$, $b_i = \mathbf{1}\{\rho(t, u_i(t)-) < 0 < \rho(t, u_i(t)+)\} - \mathbf{1}\{\rho(t, u_i(t)+) < 0 < \rho(t, u_i(t)-)\}$ and the first equation should be understood throughout $[0, T] \times \mathbb{R}$, except on the curves $u_i(t)$.

Denote by $C_0^{1,2}([0, T) \times \mathbb{R})$ the set of functions $G \colon [0, T) \times \mathbb{R} \to \mathbb{R}$ with compact support which are continuously differentiable with respect to the first variable and twice continuously differentiable with respect to the second variable. A bounded measurable function $\rho \colon [0, T] \times \mathbb{R}$ is said to be a weak solution of the Stefan problem (2.2) if for every function $G \in C_0^{1,2}([0, T) \times \mathbb{R})$,

(2.3) $$\int_{-\infty}^{+\infty} du \int_0^T dt\, \{A(\rho(t, u))\Delta G(t, u) + \omega(\rho(t, u))\partial_t G(t, u)\}$$
$$+ \int_{-\infty}^{+\infty} du\, \omega(\rho_0(u)) G(0, u) = 0,$$

where

(2.4) $$\omega(\rho) = \begin{cases} \rho - 1, & \rho < 0, \\ \rho, & \rho > 0, \\ -1, & \rho = 0. \end{cases}$$

The proof of uniqueness of weak solutions for the Stefan equation (2.2) presented in [7], Theorem 20, page 312, for boundary-valued problems can be easily adapted to our context. Furthermore, the generalized solution is continuous according to Theorem 21 of [7].

*The initial states.* Let $\mathcal{A}$ be the subset of $\Gamma$ of all configurations $(\sigma, \eta)$ for which there exists $x$ in $\mathbb{Z}$ such that

$$\sigma(z) = \begin{cases} -1, & \text{if } z \le x, \\ 1, & \text{if } z > x. \end{cases}$$

Note that $\mathcal{A}$ is stable under the dynamics induced by the generator $\mathcal{G}$.

For every $(\sigma, \eta) \in \mathcal{A}$, let

$$b = b(\sigma) := \sup\{z : \sigma(z) = -1\}$$

be the boundary of the configuration $(\sigma, \eta)$.

In the proof of the hydrodynamic behavior of the process $(\sigma, \eta)$ we impose some conditions on the initial states. Let $\{m^N, N \ge 1\}$ be a sequence of probability measures on $\Gamma$. We assume that:

(H1) For every $N \geq 1$, $m^N$ is concentrated on configurations of the set $\mathcal{A}$ such that $b(\sigma) = 0$.

(H2) There exists a bounded measurable function $\rho_0 : \mathbb{R} \to [-1, 1]$ such that

$$\int_a^\infty \rho_0(u)\, du > 0, \qquad \int_{-\infty}^{-a} \rho_0(u)\, du < 0$$

for all $a > 0$ and such that for each continuous function $G : \mathbb{R} \to \mathbb{R}$ with compact support and each $\delta > 0$,

$$\lim_{N \to \infty} m^N \left[ \left| N^{-1} \sum_{x \in \mathbb{Z}} G(x/N)\sigma(x)\eta(x) - \int_{\mathbb{R}} du\, G(u)\rho_0(u) \right| \geq \delta \right] = 0.$$

Notice that condition (H1) forces the initial profile to be negative on $\mathbb{R}_-$ and positive on $\mathbb{R}_+$.

In the case $a_{-1} = 0$, we impose one more condition:

(H3) For every $N \geq 1$,

$$m^N \{ (\sigma, \eta) \in \mathcal{A} : \eta(x) = 0, x \leq 0 \} = 1.$$

*The hydrodynamic behavior.* For each probability measure $m$ on $\Gamma$ concentrated on $\mathcal{A}$, denote by $\mathbb{P}_m^{f,N}$ the probability measure on the path space $D(\mathbb{R}_+, \Gamma)$ induced by the Markov process $(\sigma_t, \eta_t)$ with generator $\mathcal{G}$ speeded up by $N^2$ and initial measure $m$.

Denote by $\mathcal{M} = \mathcal{M}(\mathbb{R})$ the space of *signed* Radon measures on $\mathbb{R}$ endowed with the vague topology. Integration of a function $G$ with respect to a measure $\pi$ in $\mathcal{M}$ is denoted by $\langle \pi, G \rangle$. To each configuration $(\eta, \sigma) \in \Gamma$ we associate the empirical measure $\pi^N = \pi^N(\eta, \sigma)$ in $\mathcal{M}$ by assigning mass $\sigma(x) N^{-1}$ to each particle:

$$\pi^N = \frac{1}{N} \sum_{x \in \mathbb{Z}} \sigma(x) \eta(x) \delta_{x/N}.$$

Let $\pi_t^N = \pi^N(\sigma_t, \eta_t)$, $b_t^N = b(\sigma_t)/N$.

THEOREM 2.1. *Fix a sequence of initial measures $\{m^N : N \geq 1\}$ satisfying assumptions* (H1), (H2) *and* (H3) *if $a_{-1} = 0$. For each $t \geq 0$, as $N \uparrow \infty$, the empirical measure $\pi_t^N$ converges in probability to an absolutely continuous measure $\pi(t, du) = \rho(t, u)\, du$, whose density $\rho(t, u)$ is the weak solution of the Stefan problem* (2.2): *For every continuous function $G : \mathbb{R} \to \mathbb{R}$ with compact support and every $\delta > 0$*

$$\lim_{N \to \infty} \mathbb{P}_{m^N}^{f,N} \left[ \left| \langle \pi_t^N, G \rangle - \int du\, G(u)\rho(t,u) \right| \geq \delta \right] = 0.$$

*Moreover, for every $\delta > 0$*

$$\lim_{N \to \infty} \mathbb{P}^{f,N}_{m^N}[|b(t) - B(t)| > \delta] = 0,$$

*where $B$ is the solution of $B(0) = 0$,*

(2.5) $\qquad \dot{B}(t) = a_{-1}(\partial_u \rho)(t, B(t)-) - a_1(\partial_u \rho)(t, B(t)+).$

**3. The fixed boundary model.** Recall the definition of the set $\mathcal{A}$ given in the previous section and of the boundary $b = b(\sigma, \eta)$ of a configuration in $\mathcal{A}$. To each configuration $(\sigma, \eta)$ in $\mathcal{A}$ let $\xi = \xi_{\sigma,\eta}$ in $\Omega = \{0,1\}^{\mathbb{Z}}$ be the configuration viewed from the boundary:

$$\xi(x) = \eta(x + b).$$

It is not difficult to check that $\xi_t$ is a Markov process with generator $\mathcal{L}$ given by

$$\mathcal{L} = a_{-1} L_1 + L_b + a_1 L_2.$$

Here $L_1$ and $L_2$ are the parts of the generator related to the motion of particles in a simple symmetric exclusion process on $\mathbb{Z}_-$ and $\mathbb{N}$ respectively:

$$L_1 = \sum_{x \leq 0} L_{x-1,x}, \qquad L_2 = \sum_{x \geq 1} L_{x,x+1},$$

where, for every local function $f : \Omega \to \mathbb{R}$ and every integer $x$,

$$(L_{x,x+1} f)(\xi) = f(\xi^{x,x+1}) - f(\xi),$$

and $\xi^{x,x+1}$ is the configuration $\xi$ with spins at $x$, $y$ interchanged defined in (2.1).

In contrast, $L_b$ is the part of the generator related to the dissipative feature of the system: For every local function $f : \Omega \to \mathbb{R}$

$$\begin{aligned}(L_b f)(\xi) &= \xi(1)[1 - \xi(0)]\{f(\tau_{-1}(\xi - \varrho_1)) - f(\xi)\} \\ &\quad + \xi(0)[1 - \xi(1)]\{f(\tau_1(\xi - \varrho_0)) - f(\xi)\} \\ &\quad + \xi(0)\xi(1)\{f(\xi - \varrho_0 - \varrho_1) - f(\xi)\},\end{aligned}$$

where $\varrho_x$ stands for the configuration with no particles but one at $x$, and $\{\tau_x : x \in \mathbb{Z}\}$ for the group of translation so that $(\tau_x \xi)(z) = \xi(z + x)$ for all $z$ in $\mathbb{Z}$.

Fix a bounded measurable function $\lambda_0 : \mathbb{R} \to \mathbb{R}$ such that $\lambda_0(u) \geq 0$, $\lambda_0(u) \leq 0$ for $u \geq 0$, $u \leq 0$, respectively. A pair $(\lambda, D)$, where $\lambda$ is bounded measurable function $\lambda : [0, T] \times \mathbb{R} \to \mathbb{R}$ *strictly positive* a.e. on $(0, \infty)$, *strictly* *negative* a.e. on $(-\infty, 0)$ and $D : [0, T] \to \mathbb{R}$ is a bounded variation continuous function vanishing at $t = 0$, is said to be a weak solution of

$$\partial_t \lambda = a\Delta\lambda + \dot{D}(t)\partial_u \lambda,$$
(3.1) $\quad \dot{D}(t) = a_{-1}(\partial_u \lambda)(t, 0-) - a_1(\partial_u \lambda)(t, 0+),$
$$\lambda(0, \cdot) = \lambda_0(\cdot),$$

in the layer $[0, T] \times \mathbb{R}_+$, if for every function $G \in C_0^{1,2}([0, T) \times \mathbb{R})$,

(3.2)
$$\int_{-\infty}^{+\infty} du \int_0^T dt \, \{A(\lambda(t, u))\Delta G(t, u) + \lambda(t, u)\partial_t G(t, u)\}$$
$$- \int_0^T dt \int_{-\infty}^{+\infty} \{\lambda(t, u)\partial_u G(t, u) - G(t, 0)\} \, dD(t)$$
$$+ \int_{-\infty}^{+\infty} du \, \lambda_0(u) G(0, u) = 0.$$

Notice that we are requiring the solution to be *strictly* positive, negative on $(0, \infty)$, $(-\infty, 0)$, respectively.

Lemma 3.1 below shows that weak solutions of (2.2) may be obtained by a simple change of variables from weak solutions of (3.1). In particular, uniqueness of weak solutions of (3.1) follows from the uniqueness for (2.2). Indeed, assume that $(\lambda_t, D_t)$ is a weak solution of (3.1). By Lemma 3.1 and by the uniqueness of weak solutions of (2.2), $\rho(t, u) = \lambda(t, u - D_t)$ is the unique weak solution of (2.2). Since the weak solution of (2.2) is continuous and since $\lambda$ is a.e. strictly positive, negative on $(0, \infty)$, $(-\infty, 0)$, respectively, for each $t \geq 0$, $\rho(t, \cdot)$ vanishes at a unique point. This determines uniquely $D_t$ and therefore $\lambda$.

LEMMA 3.1. *Let $(\lambda_t, D_t)$ be a weak solution of* (3.1) *and let $\rho(t, u) = \lambda(t, u - D_t)$. Then, $\rho$ is a weak solution of* (2.2).

PROOF. Consider a weak solution $(\lambda_t, D_t)$ of (3.1) and write $D_t$ as the difference of two continuous increasing bounded functions: $D_t = D_t^+ - D_t^-$. Let $D^{\pm,\varepsilon}$ be smooth uniform approximations of $D^{\pm}$ and set $D_t^\varepsilon = D_t^{+,\varepsilon} - D_t^{-,\varepsilon}$ so that

(3.3) $$\lim_{\varepsilon \to 0} \sup_{0 \leq t \leq T} |D_t^\varepsilon - D_t| = 0.$$

Fix a smooth function $G : [0, T] \times \mathbb{R} \to \mathbb{R}$ with compact support and vanishing at the boundary $t = T$. Let $H^\varepsilon(t, u) = G(t, u + D_t^\varepsilon)$. $H^\varepsilon$ is a smooth function with compact support. Therefore, since $(\lambda_t, D_t)$ is a weak solution of (3.1) and since $G$ vanishes at the boundary $t = T$,

$$0 = \langle \lambda_T, G(T, u + D_T^\varepsilon) \rangle = \langle \lambda_T, H_T^\varepsilon \rangle$$

$$(3.4) \qquad = \langle \lambda_0, H_0^\varepsilon \rangle + \int_0^T ds \, \langle \lambda_s, (\partial_s + a\Delta) H_s^\varepsilon \rangle$$

$$- \int_0^T \{ \langle \lambda_s, \partial_u H_s^\varepsilon \rangle - H_s^\varepsilon(0) \} \, dD_s.$$

Recall the definition of the function $\omega$ given in (2.4) and that weak solutions of (3.1) are strictly positive in $(0, \infty)$ and strictly negative in $(-\infty, 0)$. The first two terms on the right-hand side of (3.4) can be rewritten as

$$\langle \omega(\lambda_0), H_0^\varepsilon \rangle + \int_0^T ds \, \langle \lambda_s, a\Delta H_s^\varepsilon \rangle + \int_0^T ds \, \langle \omega(\lambda_s), \partial_s H_s^\varepsilon \rangle$$

$$+ \langle \mathbf{1}\{(-\infty, 0)\}, H_0^\varepsilon \rangle + \int_0^T ds \, \langle \mathbf{1}\{(-\infty, 0)\}, \partial_s H_s^\varepsilon \rangle.$$

Since $G$, and therefore $H^\varepsilon$, vanish at the boundary $t = T$, the second line of the previous expression is equal to 0. On the other hand, the first two terms of the first line are easily seen to converge to

$$(3.5) \qquad \langle \omega(\rho_0), G_0 \rangle + \int_0^T ds \, \langle A(\rho_s), \Delta G_s \rangle$$

as $\varepsilon \downarrow 0$. The last term of the first line together with the last term on the second line of (3.4) is equal to

$$\int_0^T \langle \omega(\lambda_s), (\partial_s G)(s, u + D_s^\varepsilon) \rangle \, ds + \int_0^T \langle \omega(\lambda_s), (\partial_u G)(s, u + D_s^\varepsilon) \rangle \, dD_s^\varepsilon$$

$$- \int_0^T \{ \langle \lambda_s, (\partial_u H^\varepsilon)(s, u) \rangle - H_s^\varepsilon(0) \} \, dD_s.$$

The first term converges, as $\varepsilon \downarrow 0$, to

$$(3.6) \qquad \int_0^T ds \, \langle \omega(\rho_s), \partial_s G_s \rangle,$$

while, by definition of $\omega$ and $H^\varepsilon$, the sum of the second and third terms is equal to

$$- \int_0^T G(s, D_s^\varepsilon) \, d(D_s^\varepsilon - D_s) + \int_0^T \langle \lambda_s, \partial_u H_s^\varepsilon \rangle \, d(D_s^\varepsilon - D_s).$$

It is not difficult to show from (3.3) that this expression vanishes as $\varepsilon \downarrow 0$.

It follows from (3.4), (3.5) and (3.6) that

$$\langle \omega(\rho_0), G_0 \rangle + \int_0^T ds \, \langle A(\rho_s), \Delta G_s \rangle + \int_0^T ds \, \langle \omega(\rho_s), \partial_s G_s \rangle = 0,$$

which concludes the proof of the lemma. $\square$

*Hypotheses on the initial measures.* Fix a sequence of probability measures $\{\mu^N : N \geq 1\}$ on $\Omega$. To prove the hydrodynamic behavior of the system we will assume that

(H̃1) There exists a bounded measurable initial profile $\lambda_0 : \mathbb{R} \to [-1, 1]$ such that

$$\int_a^\infty \lambda_0(u)\, du > 0, \qquad \int_{-\infty}^{-a} \lambda_0(u)\, du < 0$$

for all $a > 0$ and such that for each $\delta > 0$ and each continuous function $G : \mathbb{R} \to \mathbb{R}$ with compact support

$$\lim_{N \to \infty} \mu^N \left[ \left| \frac{1}{N} \sum_{x \in \mathbb{Z}} G(x/N) \mathbf{1}_{\pm}(x) \xi(x) - \int_{-\infty}^{+\infty} du\, G(u) \lambda_0(u) \right| \geq \delta \right] = 0,$$

where $\mathbf{1}_{\pm}(u) = -\mathbf{1}\{u \leq 0\} + \mathbf{1}\{u > 0\}$.

In the case $a_{-1} = 0$, we impose one more condition:

(H̃2) For every $N \geq 1$,

$$\mu^N \{\xi \in \Omega : \xi(x) = 0, x \leq 0\} = 1.$$

*The hydrodynamic behavior.* For each probability measure $\mu$ on $\Omega$, denote by $\mathbb{P}^N_\mu$ the probability measure on the path space $D(\mathbb{R}_+, \Omega)$ induced by the Markov process $\xi_t$ with generator $\mathcal{L}$ speeded up by $N^2$ starting from the initial measure $\mu$.

Let $D^N_+(t)$ [resp. $D^N_-(t)$] be the total number of particles on $\mathbb{N}$ (resp. $\mathbb{Z}_-$) which left the system before time $t$ divided by $N$ and let $D^N(t) = D^N_-(t) - D^N_+(t)$. Formally $D^N_+(t) = N^{-1} \sum_{x \geq 1} \{\xi_0(x) - \xi_t(x)\}$.

THEOREM 3.2. *Fix a sequence of initial measures $\{\mu^N, N > 1\}$ satisfying (H1) and (H2) if $a_{-1} = 0$. Then, for any $t \geq 0$, any continuous $G : \mathbb{R} \to \mathbb{R}$ with compact support and any $\delta > 0$*

$$\lim_{N \to \infty} \mathbb{P}^N_{\mu^N} \left[ \left| \frac{1}{N} \sum_{x \in \mathbb{Z}} \mathbf{1}_{\pm}(x/N) G(x/N) \xi_t(x) - \int du\, G(u) \lambda(t, u) \right| \geq \delta \right] = 0,$$

$$\lim_{N \to \infty} \mathbb{P}^N_{\mu^N} [|D^N(t) - D(t)| \geq \delta] = 0,$$

*where $(\lambda, D)$ is the unique weak solution of* (3.1).

**4. Proof of Theorems 2.1 and 3.2.** We first show how Theorem 2.1 can be recovered from Theorem 3.2.

PROOF OF THEOREM 2.1. Notice first that the evolution of the original process $(\sigma_t, \eta_t)$ can be derived from the one of $\xi_t$ since $b_t^N = D^N(t)$ and $\eta_t = \tau_{-b(\sigma_t)} \xi_t$.

By Theorem 3.2, for every $t > 0$, $b_t^N$ converges in probability to $D(t)$. This proves the second statement of the theorem since $D(\cdot)$ satisfies (2.5) in virtue of (3.1).

On the other hand, if $G : \mathbb{R} \to \mathbb{R}$ is a continuous function with compact support, by the previous relations between $\xi$ and $(\eta, \sigma)$,

$$\langle \pi_t^N, G \rangle = \frac{1}{N} \sum_{x \in \mathbb{Z}} \mathbf{1}_{\pm}(x/N) G(D^N(t) + x/N) \xi_t(x).$$

By Theorem 3.2, this expression converges in probability to $\int du\, G(D(t) + u) \lambda(t, u)$ and this integral is equal to $\int du\, G(u) \rho(t, u)$ by Lemma 3.1. This concludes the proof of Theorem 2.1. □

We turn now to the proof of Theorem 3.2. We present the proof in the case $a_{-1} > 0$ which is more difficult. The same arguments apply to the case $a_{-1} = 0$.

Recall that we denote by $\mathcal{M} = \mathcal{M}(\mathbb{R})$ the space of *signed* Radon measures on $\mathbb{R}$ endowed with the vague topology and that we denote integration of a function $G$ with respect to a measure $\pi$ in $\mathcal{M}$ by $\langle \pi, G \rangle$. For each $N \geq 1$ and each configuration $\xi$ of $\Omega$, let $\pi^N$ be the empirical measure associate to $\xi$ given by

$$\pi^N = \frac{1}{N} \sum_{x \in \mathbb{Z}} \mathbf{1}_{\pm}(x) \xi(x) \delta_{x/N}.$$

Note the indicator function $\mathbf{1}_{\pm}(x)$ in the definition which corresponds to considering particles on $\mathbb{Z}_-$ as having negative charge. Let $\pi_t^N = \pi^N(\xi_t)$ and recall that we are speeding up the process by $N^2$.

Recall that $D_+^N(t)$ [resp. $D_-^N(t)$] stands for the total number of particles on $\mathbb{N}$ (resp. $\mathbb{Z}_-$) which left the system before time $t$ divided by $N$ and that $D^N(t) = D_-^N(t) - D_+^N(t)$.

With this notation, Theorem 3.2 states that the sequence $\mathbb{Q}_{\mu^N}^N$ converges weakly to the probability measure concentrated on paths $(\pi_t, D(t))$ whose first coordinate is absolutely continuous $\pi(t, du) = \lambda(t, u)\, du$, the density being the solution of (3.1) (cf. [4]). The proof consists in showing tightness, that all limit points are concentrated on absolutely continuous paths which are weak solutions of (3.1) and uniqueness of weak solutions of this equation.

Uniqueness of weak solutions of (3.1) was discussed in the previous section, while tightness is proved at the end of this section. We show now that all limit points are concentrated on weak solutions.

Fix a sequence $\mu^N$ of probability measures on $\Omega$ satisfying the assumptions of the theorem. Note that all limit points of the sequence $\mathbb{Q}^N = \mathbb{Q}^N_{\mu^N}$ are concentrated on absolutely continuous measures since in the limit the $\pi$-measure of a finite interval is bounded by its Lebesgue measure.

PROPOSITION 4.1. *All limit points $\mathbb{Q}^*$ of the sequence $\mathbb{Q}^N_{\mu^N}$ are concentrated on trajectories $(\pi_t, D_t)$ such that*

$$\langle \pi_t, G \rangle - \langle \pi_0, G \rangle = \int_0^t ds\, \langle \pi_s, (\partial_s + a\Delta)G_s \rangle - \int_0^t \{\langle \pi_s, \partial_u G_s \rangle - G(s,0)\} dD_s$$

*for every $t > 0$ and $G$ in $C_0^{1,2}([0,T) \times \mathbb{R})$.*

PROOF. The proof of this proposition is divided in several steps. We start examining some martingales associated to the empirical measure. Fix $G \in C_0^{1,2}([0,T) \times \mathbb{R})$ and $\delta > 0$. Consider the martingale $M^{G,N}$ given by

$$M_t^{G,N} = \langle \pi_t^N, G_t \rangle - \langle \pi_0^N, G_0 \rangle - \int_0^t (\partial_s + N^2 \mathcal{L}) \langle \pi_s^N, G_s \rangle\, ds.$$

An elementary computation shows that the quadratic variation $\langle M^{G,N} \rangle_t$ of this martingale is equal to the time integral of

(4.1)
$$\frac{1}{N^2} \sum_{x \neq 0} a(x/N)[\xi(x+1) - \xi(x)]^2 \{(\nabla_N G)(x/N)\}^2$$
$$+ \xi(0)[1 - \xi(1)]\{G(-1/N) + \langle \pi^N, (\nabla_N G)(x - 1/N)\rangle\}^2$$
$$+ \xi(1)[1 - \xi(0)]\{G(2/N) - \langle \pi^N, \nabla_N G \rangle\}^2$$
$$+ \xi(0)\xi(1)\{G(1/N) - G(0)\}^2.$$

In particular, since $G$ is a smooth function with compact support, by Chebyshev and Doob inequalities,

$$\mathbb{P}^N_{\mu^N}\left[\sup_{0 \leq t \leq T} |M_t^{G,N}| \geq \delta\right] \leq 4\delta^{-2} \mathbb{E}^N_{\mu^N}[\langle M^{G,N} \rangle_T],$$

$$\leq C(a,G)\delta^{-2}\left\{\frac{T}{N} + \mathbb{E}^N_{\mu^N}\left[\int_0^T \{\xi_s(0) + \xi_s(1)\}\, ds\right]\right\}$$

for some finite constant $C(a,G)$ depending only on $a(\cdot)$ and $G$. Therefore, by Lemma 5.5,

(4.2) $$\lim_{N \to \infty} \mathbb{P}^N_{\mu^N}\left[\sup_{0 \leq t \leq T} |M_t^{G,N}| \geq \delta\right] = 0.$$

On the other hand, an elementary computation shows that for every smooth function $H:\mathbb{R}\to\mathbb{R}$ with compact support,

$$\begin{aligned}
(4.3) \quad N^2\mathcal{L}\langle \pi^N, H\rangle &= \langle \pi^N, a\Delta_N H\rangle \\
&\quad + N[\xi(1)-\xi(0)]\langle \pi^N, \nabla_N H\rangle - N\xi(1)H(2/N) \\
&\quad + N\xi(0)H(-1/N) + [a_1\xi(1)+a_{-1}\xi(0)](\nabla_N H)(0) \\
&\quad + \xi(0)\langle \pi^N, \Delta_N H\rangle \\
&\quad - \xi(0)\xi(1)\{\langle \pi^N, \Delta_N H\rangle - N^{-1}(\Delta_N H)(1/N)\},
\end{aligned}$$

where $\Delta_N$ and $\nabla_N$ denote respectively the discrete Laplacian and gradient. Therefore, in view of Lemma 5.5, up to negligible terms, the martingale $M_t^{G,N}$ can be written as

$$\langle \pi_t^N, G\rangle - \langle \pi_0^N, G\rangle - \int_0^t ds\,\langle \pi_s^N, \partial_s G_s + a\Delta G_s\rangle$$
$$- \int_0^t ds\, N[\xi_s(1)-\xi_s(0)]\{\langle \pi_s^N, \partial_u G_s\rangle - G(s,0)\}.$$

By Lemma 4.3, provided we let $\varepsilon\downarrow 0$ after $N\uparrow\infty$, we may replace $N[\xi_s(0)-\xi_s(1)]$ by $\varepsilon^{-1}[D^N(s+\varepsilon)-D^N(s)]$. Therefore, in view of (4.2),

$$\lim_{\varepsilon\to 0}\limsup_{N\to\infty}\mathbb{P}_{\mu^N}^N\bigg[\sup_{0\le t\le T}\bigg|\langle \pi_t^N, G\rangle - \langle \pi_0^N, G\rangle$$
$$- \int_0^t ds\,\langle \pi_s^N, (\partial_s + a\Delta)G_s\rangle$$
$$+ \int_0^t ds\,\frac{D^N(s+\varepsilon)-D^N(s)}{\varepsilon}$$
$$\times \{\langle \pi_s^N, \partial_u G_s\rangle - G(s,0)\}\bigg| > \delta\bigg] = 0.$$

Hence, for any limit point $\mathbb{Q}^*$ of the sequence $\mathbb{Q}^N$,

$$\lim_{\varepsilon\to 0}\mathbb{Q}^*\bigg[\sup_{0\le t\le T}\bigg|\langle \pi_t, G\rangle - \langle \pi_0, G\rangle - \int_0^t ds\,\langle \pi_s, (\partial_s + a\Delta)G_s\rangle$$
$$+ \int_0^t ds\,\frac{D(s+\varepsilon)-D(s)}{\varepsilon}\{\langle \pi_s, \partial_u G_s\rangle - G(s,0)\}\bigg| > \delta\bigg] = 0.$$

By the proof of tightness of the second marginal of $\mathbb{Q}^N$ presented at the end of this section, $\mathbb{Q}^*$ is concentrated on paths $D$ which are of bounded variation. In particular, if $d^\varepsilon(s) = \varepsilon^{-1}\{D(s+\varepsilon)-D(s)\}$, for any continuous function $H:[0,T]\to\mathbb{R}$, $\int_0^t ds\, d^\varepsilon(s)H(s)$ converges, as $\varepsilon\downarrow 0$, to $\int_0^t H(s)\,dD(s)$.

Since by Lemma 4.6 $\mathbb{Q}^*$ is concentrated on paths $\pi_t$ which are continuous for the vague topology, the proposition is proved. □

The second main result of this section states that limit points of the sequence $\mathbb{Q}^N$ are concentrated on trajectories $\pi_t$ whose density is bounded below by a strictly positive function.

PROPOSITION 4.2. *For each $\delta > 0$, there exists a strictly positive continuous function $R_\delta : \mathbb{R} \to (0, 1]$ with the following property. All limit points $\mathbb{Q}^*$ of the sequence $\mathbb{Q}^N$ are concentrated on trajectories $(\pi_t, D_t)$ such that*

$$\pi_t(I) \geq \int_I R_\delta(u)\, du$$

*for all $0 \leq t \leq T$ and all finite intervals $I = [c, d]$ such that $c \geq \delta$. A similar statement holds in $(-\infty, 0)$.*

The proof of this proposition is postponed to Section 5.

LEMMA 4.3. *Fix a smooth function $G$ in $C_0^{1,2}([0, T] \times \mathbb{R})$ and $\delta > 0$. Then,*

$$\lim_{\varepsilon \to 0} \limsup_{N \to \infty} \mathbb{P}_{\mu^N}^N \left[ \sup_{0 \leq t \leq T} \left| \int_0^t ds \left\{ N \xi_s(1) - \frac{D_+^N(s+\varepsilon) - D_+^N(s)}{\varepsilon} \right\} \langle \pi_s^N, G_s \rangle \right| > \delta \right] = 0.$$

*A similar statement holds if we replace $\xi_s(1)$, $D_+^N(s)$ by $\xi_s(0)$, $D_-^N(s)$, respectively, or $\langle \pi_s^N, G_s \rangle$ by $G(s, 0)$.*

PROOF. A simple computation shows that

$$M_+^N(t) = D_+^N(t) - N \int_0^t \xi_s(1)\, ds$$

is a martingale with quadratic variation $\langle M_+^N \rangle_t$ given by $\int_0^t \xi_s(1)\, ds$. We may therefore write $\varepsilon^{-1}\{D_+^N(s+\varepsilon) - D_+^N(s)\}$ as

$$\frac{M_+^N(s+\varepsilon) - M_+^N(s)}{\varepsilon} + N \int_s^{s+\varepsilon} \xi_r(1)\, dr.$$

The martingale part is easy to estimate because it vanishes in the limit $N \uparrow \infty$. Indeed, by Chebyshev and Schwarz inequalities and by the explicit formula for the quadratic variation of the martingale $M_+^N(t)$,

$$\mathbb{Q}_{\mu^N}^N \left[ \sup_{0 \leq t \leq T} \left| \int_0^t ds\, \frac{M_+^N(s+\varepsilon) - M_+^N(s)}{\varepsilon} \langle \pi_s^N, G_s \rangle \right| > \delta \right]$$

$$\leq \frac{C(G)}{\delta \varepsilon} \mathbb{E}_{\mu^N}^N \left[ \int_0^{T+\varepsilon} ds\, |M_+^N(s)| \right] \leq \frac{C(G, T)}{\delta \varepsilon} \mathbb{E}_{\mu^N}^N \left[ \int_0^{T+\varepsilon} ds\, \xi_s(1) \right]^{1/2}.$$

Here and below $C(G)$, $C(G,T)$ are finite constants depending only on $G$ and $G$, $T$, respectively. By Lemma 5.5, this expression vanishes as $N \uparrow \infty$.

It remains to consider the difference $N\xi_s(1) - N\int_s^{s+\varepsilon} \xi_r(1)\,dr$, which is slightly more demanding. We first perform a time integration by parts to obtain that

$$\left| \int_0^t ds \left\{ N\xi_s(1) - N\int_s^{s+\varepsilon} \xi_r(1)\,dr \right\} \langle \pi_s^N, G_s \rangle \right|$$
$$\leq C(G)\left\{ N\int_0^\varepsilon ds\,\xi_s(1) + N\int_t^{t+\varepsilon} ds\,\xi_s(1) \right\}$$
$$+ \left| \int_\varepsilon^t ds\,N\xi_s(1)\left\{ \langle \pi_s^N, G_s \rangle - \frac{1}{\varepsilon}\int_{s-\varepsilon}^s \langle \pi_r^N, G_r \rangle\,dr \right\} \right|.$$

Lemma 5.6 permits to estimate the first term on the right-hand side. To estimate the second one, write the difference $\langle \pi_s^N, G_s \rangle - \langle \pi_r^N, G_r \rangle$ as

$$M_s^{G,N} - M_r^{G,N} + \int_r^s (\partial_v + N^2 \mathcal{L})\langle \pi_v^N, G_v \rangle\,dv.$$

On the one hand,

(4.4)
$$\mathbb{P}_{\mu^N}^N\left[ \sup_{0 \leq t \leq T} \left| \int_\varepsilon^t ds\,N\xi_s(1)\frac{1}{\varepsilon}\int_{s-\varepsilon}^s dr\,\{M_s^{G,N} - M_r^{G,N}\} \right| > \delta \right]$$
$$\leq \mathbb{P}_{\mu^N}^N\left[ \sup_{0 \leq t \leq T+\varepsilon} |M_t^{G,N}| \int_0^T N\xi_s(1)\,ds > \delta/2 \right].$$

Fix $\gamma > 0$. By Lemma 5.5, there exists a finite constant $A$ for which

$$\mathbb{P}_{\mu^N}^N\left[ \int_0^T N\xi_s(1)\,ds > A \right] \leq \gamma$$

for all $N \geq 1$. Therefore, the previous probability is less than or equal to

$$\gamma + \mathbb{P}_{\mu^N}^N\left[ \sup_{0 \leq t \leq T+\varepsilon} |M_t^{G,N}| > \delta/2A \right].$$

It follows from (4.2) that this expression is bounded by $\gamma$ as $N \uparrow \infty$. This proves that (4.4) vanishes in the limit $N \uparrow \infty$.

It remains to show that

(4.5) $\mathbb{P}_{\mu^N}^N\left[ \sup_{0 \leq t \leq T} \left| \int_\varepsilon^t ds\,N\xi_s(1)\frac{1}{\varepsilon}\int_{s-\varepsilon}^s dr \int_r^s (\partial_v + N^2 \mathcal{L})\langle \pi_v^N, G_v \rangle\,dv \right| > \delta \right]$

vanishes as $N \uparrow \infty$, $\varepsilon \downarrow 0$.

By the explicit expression for $(\partial_s + N^2\mathcal{L})\langle \pi_s, G_s \rangle$ given in (4.3), we have that the absolute value of the integral in this formula is dominated by

$$C(a,G) \int_0^T N\xi_s(1)\,ds\left\{ \varepsilon + \sup_{0 \leq t \leq T} \int_t^{t+\varepsilon} N[\xi_s(0) + \xi_s(1)]\,ds \right\}.$$

Repeating the argument presented just after (4.4) to eliminate $\int_0^T N\xi_s(1)\,ds$ and applying Lemma 5.6 to estimate the second term, we show that (4.5) vanishes as $N\uparrow\infty$, $\varepsilon\downarrow 0$. This concludes the proof of the lemma. □

We conclude this section by proving that the sequence of probability measures $\mathbb{Q}^N$ is tight, which in our context reduces to showing that the marginal of $\mathbb{Q}^N$ on each coordinate is tight. We start with the empirical measure. Denote by $\mathbb{Q}_1^N$ the marginal of $\mathbb{Q}^N$ on the first coordinate.

Recall that $\mathbb{Q}_1^N$ is tight if for each smooth function with compact support $G:\mathbb{R}\to\mathbb{R}$, $\langle \pi_t^N, G\rangle$ is tight as a random sequence on $D(\mathbb{R}_+,\mathbb{R})$. Now fix such a function. To prove tightness for $\langle \pi_t^N, G\rangle$ it is enough to verify the following two conditions:

(i) The finite-dimensional distributions of $\langle \pi_t^N, G\rangle$ are tight.
(ii) For every $\delta > 0$
$$\lim_{\varepsilon\to 0}\limsup_{N\to\infty}\mathbb{P}_{\mu^N}^N\left[\sup_{|s-t|\le\varepsilon}|\langle \pi_t, G\rangle - \langle \pi_s, G\rangle| > \delta\right] = 0.$$

Condition (i) is a trivial consequence of the fact that the empirical measure has finite total mass on any compact interval. In order to prove condition (ii), consider the martingale with respect to $\mathcal{F}$ given by

$$M_t^{G,N} = \langle \pi_t^N, G\rangle - \langle \pi_0^N, G\rangle - \int_0^t N^2 \mathcal{L}\langle \pi_s^N, G\rangle\,ds.$$

Here the index $N$ indicates that we are considering the process speeded up by $N^2$. Therefore,

$$\langle \pi_t^N, G\rangle - \langle \pi_s^N, G\rangle = M_t^{G,N} - M_s^{G,N} + \int_s^t N^2 \mathcal{L}\langle \pi_r^N, G\rangle\,dr.$$

From the previous expression, condition (ii) is a consequence of the next two lemmas.

LEMMA 4.4. *For every $\delta > 0$ and every function $G$ in $C^2(\mathbb{R})$ with compact support,*

$$\lim_{\varepsilon\to 0}\limsup_{N\to\infty}\mathbb{P}_{\mu^N}^N\left[\sup_{|s-t|\le\varepsilon}\left|\int_s^t N^2\mathcal{L}\langle \pi_r^N, G\rangle\,dr\right| > \delta\right] = 0.$$

PROOF. In view of (4.3), since $G$ is a smooth function with compact support, the expression inside the absolute value is bounded above by

$$C(a,G)\left\{\varepsilon + \int_s^t N\{\xi_r(0) + \xi_r(1)\}\,dr\right\}$$

for some finite constant which depends on $a$, $G$ only. To conclude the proof, it remains to recall Lemma 5.6. □

LEMMA 4.5. *For every function $G$ in $C^2(\mathbb{R})$ with compact support and every $\delta > 0$,*

$$\lim_{\varepsilon \to 0} \limsup_{N \to \infty} \mathbb{P}^N_{\mu^N}\left[\sup_{|s-t|\leq\varepsilon} |M^{G,N}_t - M^{G,N}_s| > \delta\right] = 0.$$

PROOF. Denote by $\langle M^{G,N}\rangle_t$ the quadratic variation of the martingale $M^{G,N}_t$. By the Doob inequality,

$$\mathbb{P}^N_{\mu^N}\left[\sup_{|s-t|\leq\varepsilon} |M^{G,N}_t - M^{G,N}_s| > \delta\right] \leq \mathbb{P}^N_{\mu^N}\left[\sup_{0\leq t\leq T} |M^{G,N}_t| > \delta/2\right]$$

$$\leq \frac{4}{\delta^2}\mathbb{E}^N_{\mu^N}[\langle M^{G,N}\rangle_T].$$

By the explicit expression for $\langle M^{G,N}\rangle_T$ given in (4.1), the previous expression is bounded by

$$\frac{C(a,G)}{\delta^2}\left\{\frac{T}{N} + \mathbb{E}^N_{\mu^N}\left[\int_0^T \{\xi_s(0) + \xi_s(1)\}\,ds\right]\right\}.$$

To conclude the proof of the lemma, it remains to apply Lemma 5.5. □

We turn now to the tightness of the second marginal of $\mathbb{Q}^N$. Since $D^N_-(0) = D^N_+(0) = 0$, we need only to show that

$$\limsup_{\varepsilon\to 0}\limsup_{N\to\infty} \mathbb{P}^N_{\mu^N}\left[\sup_{|s-t|\leq\varepsilon} |D^N_+(t) - D^N_+(s)| > \delta\right] = 0$$

for every $\delta > 0$ and a similar statement for $D^N_-(t)$ in place of $D^N_+(t)$. In fact, we claim that for every $\delta > 0$, there exists $\varepsilon > 0$ such that

(4.6) $$\limsup_{N\to\infty} \mathbb{P}^N_{\mu^N}\left[\sup_{|s-t|\leq\varepsilon} |D^N_+(t) - D^N_+(s)| > \delta\right] = 0.$$

Since $D^N_+(t)$ is increasing, the previous probability is bounded above by

$$\sum_{j=1}^{T\varepsilon^{-1}} \mathbb{P}^N_{\mu^N}[D^N_+([j+1]\varepsilon) - D^N_+(j\varepsilon) > \delta/2].$$

It follows from Lemma 5.4 of the next section that for each $\delta > 0$, there exists $\varepsilon > 0$ such that

$$\limsup_{N\to\infty} \mathbb{P}^N_{\mu^N}[D^N_+(t+\varepsilon) - D^N_+(t) > \delta] = 0$$

uniformly in $0 \leq t \leq T$. This proves that the second marginal of $\mathbb{Q}^N$ is tight.

We summarize in the next lemma what we just obtained. Notice that we proved tightness in the uniform topology.

LEMMA 4.6. *The sequence $\mathbb{Q}^N$ is tight in the uniform topology. In particular, all limit points are concentrated on continuous trajectories for the vague topology.*

PROOF OF THEOREM 3.2. By Lemma 4.6, the sequence is tight in the uniform topology and by Propositions 4.1 and 4.2, all limit points are concentrated on weak solutions $(\lambda, D)$ of (3.1). By uniqueness of weak solutions, presented in Section 3, the theorem is proved. □

**5. Coupling.** We prove in this section some important but technical results which are used in the previous sections. The proofs rely on a coupling between the process $\xi_t$ defined in Section 3 and an exclusion process $\zeta_t$ similar to $\xi_t$ with the difference that the configuration is not translated when a single particle dies at the boundary. The generator $\mathcal{L}'$ of this process is therefore $a_{-1}L_{-1} + a_1 L_1 + L'_b$, where $L'_b$ is given by

$$\begin{aligned}(L'_b f)(\zeta) &= \zeta(1)[1 - \zeta(0)]\{f(\zeta - \varrho_1) - f(\zeta)\}\\ &\quad + \zeta(0)[1 - \zeta(1)]\{f(\zeta - \varrho_0) - f(\zeta)\}\\ &\quad + \zeta(0)\zeta(1)\{f(\zeta - \varrho_0 - \varrho_1) - f(\zeta)\}.\end{aligned}$$

Notice that both marginal processes on $\mathbb{Z}_-$ and on $\mathbb{N}$ behave as an exclusion process with disappearance at the boundary, whose hydrodynamic behavior is well known. The leading idea of this section is to show, through appropriate couplings, that the original process does not differ much in several aspects from the one defined above.

Denote by $\{\zeta_t : t \geq 0\}$ the Markov process with generator $\mathcal{L}'$ speeded up by $N^2$ and recall that we denote by $\{\xi_t : t \geq 0\}$ the Markov process with generator $\mathcal{L}$ speeded up by $N^2$. Let $D^+_\xi([s,t])$ be the total number of $\xi$-particles in $\mathbb{N}$ which died in the time interval $[s,t]$. $D^-_\xi([s,t])$, $D^\pm_\zeta([s,t])$ are defined analogously.

LEMMA 5.1. *There exists a coupling $(\xi_t, \zeta_t)$ for which*

(5.1) $\qquad D^-_\xi([0,t]) + D^+_\xi([0,t]) \leq 2D^-_\zeta([0,t]) + 2D^+_\zeta([0,t])$

*for all $t \geq 0$.*

PROOF. The coupling between $\xi_t$ and $\zeta_t$ can be described as follows. Assume that the initial configurations are identical at time 0: $\zeta_0 = \xi_0$. Label all particles and denote by $X^j_t$ (resp. $Y^j_t$) the position at time $t$ of the $j$th $\xi$- (resp. $\zeta$-) particle. We assume that $X^j_0 < X^k_0$ if $j < k$, $X^0_0 \leq 0 < X^1_0$, $X^j_0 = Y^j_0$ for all $j$.

The $\xi$- and $\zeta$-particles with the same label jump together preserving the order of the labels until a particle dies. If two $\xi$-particles die simultaneously, we couple the disappearance of the $\xi$- and $\zeta$-particles. If it is a single $\xi$-particle which disappears, assume, without loss of generality, that it has a positive label and denote this time by $T_0$. Note that due to the translation, $X_{T_0}^j = Y_{T_0}^j + 1$ for all labels $j$ associated to alive particles.

Denote by $T_1$ the first time after $T_0$ in which the total number of disappearances of $\xi$-particles in $\mathbb{N}$ is equal to the total number of disappearances of $\xi$-particles in $\mathbb{Z}_-$:

$$T_1 = \inf\{t > T_0 : D_\xi^+([T_0,t]) = D_\xi^-([T_0,t])\}.$$

In the time interval $[T_0,T_1]$, the coupling forces the $\xi$- and $\zeta$-particles with the same label to jump together. It may happen, however, that a $\zeta$-particle on $\mathbb{N}$ disappears while its corresponding $\xi$-particle remains alive. In this case, the $\xi$-particle becomes a second-class particle to allow the coupled particles to jump together. The same phenomenon may occur on $\mathbb{Z}_-$, where a $\xi$-particle may disappear while its corresponding $\zeta$-particle remains alive. In this case also the $\zeta$-particle becomes a second-class particle. Notice, in particular, that the difference $X_t^j - Y_t^j$ does not depend on $j$ for coupled particles.

Due to the translations, for any $T_0 \leq t \leq T_1$, the total number of $\xi$-particles which died on $\mathbb{N}$ is bounded by the total number of $\zeta$-particles which died on $\mathbb{N}$:

$$D_\xi^+([T_0,t]) \leq D_\zeta^+([T_0,t]).$$

Also due to our definition of $T_1$, for any $T_0 \leq t \leq T_1$, the total number of $\xi$-particles which died on $\mathbb{Z}_-$ is bounded by the total number of $\xi$-particles which died on $\mathbb{N}$:

$$D_\xi^-([T_0,t]) \leq D_\xi^+([T_0,t]).$$

Therefore, (5.1) holds for any $0 \leq t \leq T_1$.

To proceed by iteration, let $K_1 = D_\zeta^+([T_0,T_1])$ be the total number of $\zeta$-particles which died in $\mathbb{N}$ in the time interval $[T_0,T_1]$, and let $k_1 = D_\xi^+([T_0,T_1])$ be the total number of $\xi$-particles which died in $\mathbb{N}$ in the time interval $[T_0,T_1]$. By definition of the coupling, at time $T_1$, there are $K_1 - k_1$ uncoupled $\xi$-particles on $\mathbb{N}$ and $K_1 - k_1$ $\zeta$-particles which died on $\mathbb{N}$ associated to the $K_1 - k_1$ uncoupled $\xi$-particles. These $K_1 - k_1$ $\zeta$-particles should not be forgotten, since they will be used to compensate the eventual death of the second-class uncoupled $\xi$-particles. The important fact for the recurrence argument is that the number of dead $\zeta$-particles which were not used to compensate the death of $\xi$-particles is at least equal to the number of second-class uncoupled $\xi$-particles. Notice also that there might be second-class uncoupled $\zeta$-particles on $\mathbb{Z}_-$. Since they do not play any role in the argument, we do not refer to them again.

Assume that the first single $\xi$-particle to die after $T_1$ is in $\mathbb{Z}_-$. If it were in $\mathbb{N}$, we could repeat the arguments presented in the last paragraph to arrive at the same conclusions obtained there and iterate again the argument.

Denote by $T_2$ the first time after $T_1$ in which the total number of disappearances of $\xi$-particles in $\mathbb{N}$ is equal to the total number of disappearances in $\mathbb{Z}_-$:

$$T_2 = \inf\{t > T_1 : D_\xi^+([T_1,t]) = D_\xi^-([T_1,t])\}.$$

The coupling in $[T_1, T_2]$ is the same described before, in which coupled particles jump together until one of them dies. For $T_1 < t \leq T_2$, let $L_2(t) = D_\zeta^-([T_1,t])$, $L_2 = L_2(T_2)$, $\ell_2(t) = D_\xi^-([T_1,t])$, $\ell_2 = \ell_2(T_2)$. By definition of the coupling $D_\xi^+([T_1,t]) \leq \ell_2(t) \leq L_2(t)$, so that (5.1) holds for $0 \leq t \leq T_2$.

On the other hand, at time $T_2$, there are:

(a) at most $K_1 - k_1$ uncoupled $\xi$-particles on $\mathbb{N}$. There might be less since a second-class $\xi$-particle might have died, its death being compensated by the death of a $\zeta$-particle on $\mathbb{Z}_-$,

(b) $L_2 - \ell_2$ uncoupled $\xi$-particles on $\mathbb{Z}_-$,

(c) $K_1 - k_1$ $\zeta$-particles which died on $\mathbb{N}$ (in the time interval $[T_0, T_1]$) and which are associated to the remaining uncoupled $\xi$-particles on $\mathbb{N}$,

(d) $L_2 - \ell_2$ $\zeta$-particles which died on $\mathbb{Z}_-$ associated to the remaining uncoupled $\xi$-particles on $\mathbb{Z}_-$.

Thus, at time $T_2$, the total number of uncoupled second-class $\xi$-particles is still smaller than the total number of dead and disassociated $\zeta$-particles.

Assume, without loss of generality, that the first single particle to die after $T_2$ is in $\mathbb{N}$. Denote by $T_3$ the first time after $T_2$ in which the total number of disappearances of $\xi$-particles on $\mathbb{N}$ is equal to the total number of disappearances in $\mathbb{Z}_-$:

$$T_3 = \inf\{t > T_2 : D_\xi^+([T_2,t]) = D_\xi^-([T_2,t])\}.$$

The coupling remains the same. The next argument, though elementary, requires much notation. For $T_2 < t \leq T_3$, let:

(a) $j_3(t)$ be the total number of second-class $\xi$-particles which die on $\mathbb{N}$ in the time interval $[T_2, t]$,

(b) $k_3(t) - j_3(t)$ be the total number of first-class $\xi$-particles which die on $\mathbb{N}$ in the time interval $[T_2, t]$,

(c) $K_3(t) - j_3(t)$ be the total number of $\zeta$-particles which die on $\mathbb{N}$ in the time interval $[T_2, t]$,

(d) $\ell_3(t)$ be the total number of $\xi$-particles which die on $\mathbb{Z}_-$ in the time interval $[T_2, t]$.

By the definition of the coupling, $\ell_3(t) \leq k_3(t) \leq K_3(t)$. The $j_3(t)$ second-class $\xi$-particles which died on $\mathbb{N}$ were associated to $\zeta$-particles which died before. Since there is a factor 2 in (5.1), we also associate to these $\zeta$-particles, $j_3(t) \wedge \ell_3(t)$ $\xi$-particles which died on $\mathbb{Z}_-$. The $k_3(t) - j_3(t)$ first-class $\xi$-particles which died on $\mathbb{N}$ are taken care of by $k_3(t) - j_3(t)$ $\zeta$-particles which died on $\mathbb{N}$. The factor 2 in (5.1) allows to include $(\ell_3(t) - j_3(t))^+ \leq k_3(t) - j_3(t)$ $\xi$-particles which died on $\mathbb{Z}_-$. Up to this point we showed that all disappearances of $\xi$-particles in $[T_2, t]$ can be compensated by disappearance of $\zeta$-particles in $[T_2, t]$ and by disassociated $\zeta$-particles which died before $T_2$. Therefore, (5.1) holds in the time interval $[T_2, T_3]$.

To be able to iterate this argument notice that there are $K_3(T_3) - k_3(T_3)$ second-class $\xi$-particles created on $\mathbb{N}$ in the time interval $[T_2, T_3]$ and $K_3(T_3) - k_3(T_3)$ $\zeta$-particles which died in this interval and whose deaths were not used to compensate $\xi$ deaths. We may therefore associate these new second-class $\xi$-particles to these newly dead $\zeta$-particles and iterate the argument. This concludes the proof of the lemma. $\square$

Denote by $\zeta_t^-$, $\zeta_t^+$ the marginals of the process $\zeta_t$ on $\mathbb{Z}_-$, $\mathbb{N}$, respectively. Notice that both marginals evolve as an exclusion process in which particles leave the system at the boundary. This system plays an important role in the sequel and deserves a notation. For $b > 0$, denote by $\beta_t$ the Markov process on $\{0,1\}^{\mathbb{Z}_+}$ with generator $\mathfrak{L} = \mathfrak{L}_b$ given by

$$(\mathfrak{L}f)(\beta) = b \sum_{x \geq 0} \{f(\beta^{x,x+1}) - f(\beta)\} + \beta(0)\{f(\beta - \varrho_0) - f(\beta)\}.$$

For $T > 0$ and a measure $\mu$ on $\{0,1\}^{\mathbb{Z}_+}$, denote by $\tilde{\mathbb{P}}_\mu^N$ the probability on $D([0,T], \{0,1\}^{\mathbb{Z}_+})$ induced by the Markov process $\beta_t$ speeded up by $N^2$ starting from $\mu$. Expectation with respect to $\tilde{\mathbb{P}}_\mu^N$ is denoted by $\tilde{\mathbb{E}}_\mu^N$.

It is well known that the process $\beta_t$ has a hydrodynamic description. Let $\tilde{\pi}_t^N$ be the empirical measure associated to $\beta_t$: $\tilde{\pi}_t^N = N^{-1} \sum_{x \geq 0} \beta_t(x) \delta_{x/N}$.

PROPOSITION 5.2. *Consider a sequence of probability measures $\tilde{\mu}^N$ on $\{0,1\}^{\mathbb{Z}_+}$ such that*

$$\lim_{N \to \infty} \tilde{\mu}^N \left[ \left| \langle \pi^N, G \rangle - \int_0^\infty \rho_0(u) G(u) \, du \right| > \delta \right] = 0$$

*for some measurable function $\rho_0 : \mathbb{R}_+ \to [0,1]$, every $\delta > 0$ and every continuous function with compact support $G$. Then,*

$$\lim_{N \to \infty} \tilde{\mathbb{P}}_{\tilde{\mu}^N} \left[ \left| \langle \pi_t^N, G \rangle - \int_0^\infty \rho(t,u) G(u) \, du \right| > \delta \right] = 0$$

*for every $\delta > 0$ and continuous function with compact support $G$, where $\rho$ is the solution of the linear equation*

(5.2)
$$\partial_t \rho = b\Delta\rho \quad on\ \mathbb{R}_+,$$
$$\rho(t,0) = 0, \qquad \rho(0,\cdot) = \rho_0(\cdot).$$

The proof of this result is similar to the one of Theorem 2.1 in [5]. Moreover, the solution of (5.2) can be represented in terms of a standard Brownian motion $W_t$ with absorption at the boundary $u = 0$:

(5.3)
$$\rho(t,u) = E_u[\rho_0(\sqrt{b}W_t)].$$

The coupling presented in Lemma 5.1 together with the hydrodynamic behavior stated in Proposition 5.2 permit to estimate the total number of particles which left the system in the original process $\xi_t$. This is the content of the next two lemmas. Recall the definition of $D_+^N(t)$, $D_-^N(t)$ introduced in Section 4. Denote by **1** the probability measure on $\{0,1\}^{\mathbb{Z}_+}$ concentrated on the configuration with all sites occupied and by $D_\beta(t)$ the total number of $\beta$-particles which left the system before time $t$.

LEMMA 5.3. *There exists a finite constant $C_0$ depending only on $a_1$, $a_{-1}$ such that*

$$\limsup_{N\to\infty} \mathbb{E}_\mu^N[D_+^N(t)] \leq C_0\sqrt{t}$$

*for all $t \geq 0$ and all probability measures $\mu$. The lemma remains in force if $D_+^N$ is replaced by $D_-^N$.*

PROOF. By Lemma 5.1, the expectation in the statement of the lemma is bounded above by $2\mathbb{E}_\mu^N[D_\zeta(t)]$, where $D_\zeta(t)$ stands for the total number of particles which left the system in the time interval $[0,t]$ for the $\zeta$ process. Since both marginals of $\zeta$ evolve as an exclusion process with disappearance at the boundary, the previous expectation is bounded above by $4\max_{b=a_1,a_{-1}} \hat{\mathbb{E}}_\mu^N[D_\beta^N(t)]$. By monotonicity, this latter expectation is less than or equal to $4\max_{b=a_1,a_{-1}} \tilde{\mathbb{E}}_\mathbf{1}^N[D_\beta^N(t)]$. By the hydrodynamic limit of $\beta$, this expectation converges, as $N \uparrow \infty$, to

$$4\max_{b=a_1,a_{-1}} \int_0^\infty \{1 - \rho_b(t,u)\}\,du,$$

where $\rho_b$ is the solution of (5.2) with initial condition $\rho_0$ constant equal to 1. With this initial condition, the solution of this equation can be written as $\rho_b(t,u) = 1 - 2P[B_t \geq u/\sqrt{2b}]$, where $B_t$ is a standard Brownian motion. In particular, the previous displayed equation is equal to

$$4\sqrt{2t}\max\{\sqrt{a_{-1}}, \sqrt{a_1}\}E[|B_1|].$$

This concludes the proof of the lemma. □

LEMMA 5.4. *For every $\delta > 0$, $T > 0$, there exists $\varepsilon > 0$ such that*

$$\limsup_{N \to \infty} \mathbb{P}^N_{\mu^N}[D^N_+(t+\varepsilon) - D^N_+(t) > \delta] = 0$$

*uniformly for $0 \leq t \leq T$. The statement remains in force if we replace $D^N_+$ by $D^N_-$.*

PROOF. Fix $0 \leq t \leq T$. Denote by $\mu^N(t)$ the state of the process at time $t$. With this notation the probability appearing in the statement can be written as

$$\mathbb{P}^N_{\mu^N(t)}[D^N_+(\varepsilon) > \delta].$$

By Lemma 5.1 and by attractiveness of the $\beta$-process, the previous expression is less than or equal to

$$2 \max_{b=a_1, a_{-1}} \tilde{\mathbb{P}}^N_{\mathbf{1}}[D^N_\beta(\varepsilon) > \delta/4].$$

By the hydrodynamic limit of $\beta$ and the proof of the previous lemma, $D^N_\beta(\varepsilon)$ converges in probability to

$$\int_0^\infty \{1 - \rho_b(\varepsilon, u)\}\, du = \sqrt{2b\varepsilon} E[|B_1|].$$

In particular, if $\varepsilon$ is chosen small enough for the last expression to be less than $\delta/4$, the previous probability vanishes as $N \uparrow \infty$. This concludes the proof of the lemma. □

LEMMA 5.5. *For every $t \geq 0$,*

$$\sup_N \mathbb{E}^N_{\mu^N}\left[\int_0^t N\{\xi_s(0) + \xi_s(1)\}\, ds\right] < \infty.$$

PROOF. Recall the definition of the martingale $M^N_+(t)$ introduced in the beginning of the proof of Lemma 4.3. In particular, $\mathbb{E}^N_{\mu^N}[\int_0^t N\xi_s(1)\, ds] = \mathbb{E}^N_{\mu^N}[D^N_+(t)]$. It remains to recall Lemma 5.3 to estimate the expectation appearing in the statement of the lemma for large $N$. For small $N$, it is enough to bound $D^N_+(t)$ by a Poisson point process. □

LEMMA 5.6. *For every $T > 0$ and $\delta > 0$,*

$$\lim_{\varepsilon \to 0} \limsup_{N \to \infty} \mathbb{P}^N_{\mu^N}\left[\sup_{0 \leq t \leq T} \int_t^{t+\varepsilon} N\{\xi_s(0) + \xi_s(1)\}\, ds > \delta\right] = 0.$$

PROOF. Recall the definition of the martingale $M_+^N(t)$ introduced in the beginning of the proof of Lemma 4.3. With this notation, in order to prove the lemma we need to show that

$$\lim_{\varepsilon \to 0} \limsup_{N \to \infty} \mathbb{P}_{\mu^N}^N \left[ \sup_{0 \leq t \leq T} \{M_+^N(t+\varepsilon) - M_+^N(t)\} > \delta \right] = 0,$$

$$\lim_{\varepsilon \to 0} \limsup_{N \to \infty} \mathbb{P}_{\mu^N}^N \left[ \sup_{0 \leq t \leq T} \{D_+^N(t+\varepsilon) - D_+^N(t)\} > \delta \right] = 0,$$

and a similar statement with $M_-^N(t)$, $D_-^N(t)$ in place of $M_+^N(t)$, $D_+^N(t)$.

The martingale part is easy. By Doob's inequality and the explicit formula for the quadratic variation of $M_+^N(t)$ given at the beginning of the proof of Lemma 4.3, the probability which needs to be estimated is less than or equal to

$$\mathbb{P}_{\mu^N}^N \left[ \sup_{0 \leq t \leq T+\varepsilon} |M_+^N(t)| > \delta/2 \right] \leq \frac{4}{\delta^2} \mathbb{E}_{\mu^N}^N \left[ \int_0^{T+\varepsilon} \xi_s(1)\,ds \right].$$

By the previous lemma, this expression vanishes as $N \uparrow \infty$ for any $\varepsilon > 0$, $\delta > 0$.

On the other hand, the jump part has been estimated just after the proof of Lemma 4.5. $\square$

We conclude this section with the proof of Proposition 4.2 which relies on the following lemma. For a subset $I$ of $\mathbb{R}$, denote by $\mathbf{1}\{I\}$ the indicator function of $I$.

LEMMA 5.7. *Fix $T > 0$ and let $\mu^N$ be a probability measure satisfying assumption ($\tilde{H}1$). There exist a finite constant $A_0$ and a strictly positive continuous function $R_T : (-\infty, -A_0] \cup [A_0, \infty) \to \mathbb{R}_+$ such that*

$$\limsup_{N \to \infty} \mathbb{P}_{\mu^N}^N \left[ |\langle \pi_t^N, \mathbf{1}\{[a,b]\} \rangle| < \int_a^b R_T(u)\,du \right] = 0$$

*for all $0 \leq t \leq T$ and all intervals $[a,b]$ such that $a > A_0$ or $b < -A_0$.*

PROOF. Fix $T > 0$, $0 \leq t \leq T$ and let $|D^N|(t) = D_+^N(t) + D_-^N(t)$ be the total number of particles which left the system before time $t$ divided by $N$. By Lemma 5.4, there exists $A_1 > 0$ such that

$$\limsup_{N \to \infty} \mathbb{P}_{\mu^N}^N [|D^N|(T) > A_1] = 0.$$

Fix such $A_1 > 0$ and couple the $\xi$ process with a $\zeta$ process as in Lemma 5.1 with the additional property that $\zeta$-particles which jump to the interval $\{-A_1 N, \ldots, A_1 N\}$ are removed. In particular, all $\xi$-particles initially in this interval become instantaneously second-class particles.

On the set $|D^N|(T) \leq A_1$, we have that a particle $Y_t^j$ is alive only if the particle $X_t^j$ is alive and that the distance $|X_t^j - Y_t^j|$, which does not depend on $j$ for alive particles, is bounded by $N|D^N|(T)$. Therefore, if we fix an interval $I = [a,b]$ with $a > 3A_1$ or $b < -3A_1$, on the set $|D^N|(T) \leq A_1$,

$$\langle \pi^N(\xi_t), \mathbf{1}\{I\} \rangle = N^{-1} \sum_{x/N \in I} \xi_t(x) = N^{-1} \sum_j \mathbf{1}\{X_t^j \in NI\},$$

where the last summation is performed over all indices $j$ corresponding to alive particles, and is bounded below by

$$\inf_{|v| \leq A_1} N^{-1} \sum_j \mathbf{1}\{Y_t^j \in N(v+I)\} = \inf_{|v| \leq A_1} \langle \pi^N(\zeta_t), \mathbf{1}\{v+I\} \rangle.$$

In these formulas, $v + I = \{u + v : u \in I\}$ and $NJ = \{Nu : u \in J\}$.

Since $\inf_{|v| \leq A_1} \langle \pi, \mathbf{1}\{v+I\} \rangle$ is a continuous function for the vague topology because all measures have density bounded by 1, by the hydrodynamic limit for $\zeta$, the last expression converges in probability to

$$\inf_{|v| \leq A_1} \langle \rho_t, \mathbf{1}\{v+I\} \rangle \geq \int_a^b \inf_{|v-u| \leq A_1} \rho(t,v) \, du,$$

where $\rho$ is the solution of the linear equation (5.2) with initial condition $\rho_0$ and boundary condition $\rho(t, \pm A_1) = 0$. By the explicit formula (5.3) for the solution of (5.2), $\rho(t,u)$ is smooth. Moreover, since $\rho_0$ satisfies condition ($\tilde{\text{H}}$1), $\rho(t,u)$ is strictly positive. Therefore, for each $u > 3A_1$,

$$R(T,u) = \inf_{0 \leq t \leq T} \inf_{|v-u| \leq A_1} \rho(t,v) > 0.$$

This concludes the proof of the lemma. □

PROOF OF PROPOSITION 4.2. In Lemma 5.7 we proved the proposition for intervals far from the origin. We estimate now the density on intervals close to the origin repeating the same argument presented in the proof of Lemma 5.7 and using the fact that the total number of particles which leave the system in a small time interval cannot be too large.

Fix $T > 0$ and recall from (4.6) that for each $\delta > 0$ there exists $\varepsilon = \varepsilon(\delta) > 0$ such that

$$\limsup_{N \to \infty} \mathbb{P}_{\mu^N}^N \left[ \sup_{0 \leq t \leq T} |D^N|(t+\varepsilon) - |D^N|(t) > \delta \right] = 0.$$

We shall estimate without loss of generality the density on $\mathbb{R}_+$. Fix $a > 0$. Let $\delta = a/3$, $s = t - \varepsilon(\delta)$ and denote by $\mu^N(s)$ the state of the $\xi$ process at time $s$. By Lemma 5.7, at time $s$, the particles' density is bounded below by a strictly positive function $R_T$:

$$\limsup_{N \to \infty} \mu^N(s)[\langle \pi^N, \mathbf{1}\{I\} \rangle < R_T(I)] = \limsup_{N \to \infty} \mathbb{P}_{\mu^N}^N[\langle \pi_s^N, \mathbf{1}\{I\} \rangle < R_T(I)] = 0$$

for all intervals $I = [c, d]$ with $c > A_0$. Here, $R_T(I) = \int_I R_T(u)\, du$.

Starting from $\mu^N(s)$, we couple the $\xi$ process with a $\zeta$ process as in the proof of Lemma 5.1 with the additional feature that $\zeta$-particles which reach the set $\{1, \ldots, \delta N\}$ are killed. Following the argument presented in the proof of Lemma 5.7, on the set $|D^N|(t) - |D^N|(s) \le \delta$ we obtain that

$$\langle \pi^N(\xi_t), \mathbf{1}\{J\}\rangle \ge \inf_{|v| \le \delta} \langle \pi^N(\zeta_t), \mathbf{1}\{v + J\}\rangle,$$

for every interval $J = [c, d]$ with $c \ge a$. The asymptotic behavior of the right-hand side of this inequality is given by the hydrodynamic limit of the $\zeta$ process in the time interval $[s, t]$. The fact that we do not know the law of large number for the empirical measure at time $s$ is not a problem. In fact, it is not difficult to show that $\pi^N(\zeta_r)$ is a tight sequence and that all limit points are concentrated on weak solutions of (5.2) with boundary condition $\rho(r, \delta) = 0$ for $s \le r \le t$ and on trajectories $\pi$ which at time $s$ are bounded below by $R_T$ on the interval $[A_0, \infty)$. By monotonicity of weak solutions of (5.2), at time $t$ the empirical measure has a density bounded below by the solution of (5.2) with initial condition $R_T \mathbf{1}\{[A_0, \infty)\}$. Therefore, in the limit $N \uparrow \infty$, the right-hand side of the previous equation is bounded below by

$$\int_c^d \inf_{|v-u| \le \delta} \rho(\varepsilon, v)\, du,$$

where $\rho$ is the solution of (5.2) with initial condition $\rho_0 = R_T \mathbf{1}\{[A_0, \infty)\}$ and boundary condition $\rho(r, \delta) = 0$ for $0 \le r \le \varepsilon$. Here again, the explicit formula (5.3) for the solution of (5.2) shows that the continuous function $\inf_{|v-u| \le \delta} \rho(\varepsilon, v)$ is strictly positive on $[a, \infty)$. This concludes the proof of the proposition. $\square$

## REFERENCES

[1] BERTINI, L., BUTTÀ, P. and RÜDIGER, B. (2000). Interface dynamics and Stefan problem from microscopic conservative model. *Rend. Mat. Appl.* **19** 547–581. MR1789487

[2] CHAYES, L. and SWINDLE, G. (1996). Hydrodynamic limits for one-dimensional particle systems with moving boundaries. *Ann. Probab.* **24** 559–598. MR1404521

[3] VALLE, G. (2003). Evolution of the interfaces in a two dimensional Potts model. Preprint.

[4] KIPNIS, C. and LANDIM, C. (1997). *Scaling Limits of Interacting Particle Systems.* Springer, Berlin. MR1707314

[5] LANDIM, C., OLLA, S. and VOLCHAN, S. (1998). Driven tracer particle and Einstein relation in one dimensional symmetric simple exclusion process. *Comm. Math. Phys.* **192** 287–307. MR1617558

[6] REZAKHANLOU, F. (1990). Hydrodynamic limit for a system with finite range interactions. *Comm. Math. Phys.* **129** 445–480. MR1051500

[7] RUBINSTEIN, L. I. (1971). *The Stefan Problem*. Amer. Math. Soc., Providence, RI. MR0351348

[8] STEFAN, J. (1889). Über einige Probleme der Theorie der Wärmeleitung. *S-B Wien Akad. Mat. Natur.* **98** 173–484.

IMPA  
ESTRADA DONA CASTORINA 110  
CEP 22460 RIO DE JANEIRO  
BRASIL  
AND  
CNRS UMR 6085  
UNIVERSITÉ DE ROUEN  
76128 MONT SAINT AIGNAN  
FRANCE  
E-MAIL: landim@impa.br

ECOLE POLYTECHNIQUE FÉDÉRALE DE LAUSANNE  
LAUSANNE  
SUISSE  
E-MAIL: glauco.valle@epfl.ch